\documentclass[12pt]{amsart}

\makeatletter \@addtoreset{equation}{section}
\@addtoreset{figure}{section}

\makeatother

\usepackage{xypic,amscd}
\usepackage{amssymb}

\newcommand{\PP}{\mathbb{P}}

\newcommand{\FF}{\mathbb{F}}

\newcommand{\CC}{\mathbb{C}}
\newcommand{\QQ}{\mathbb{Q}}

\newtheorem{theorem}[equation]{Theorem}

\newtheorem{lemma}[equation]{Lemma}

\theoremstyle{definition}

\newtheorem{definition}[equation]{Definition}

\newtheorem{remark}[equation]{Remark}
\newtheorem{zero}[equation]{}
\newtheorem{case}{}[equation]

\theoremstyle{remark}

\date{}

\begin{document}

\title[]{The maximal number of singular points on log del Pezzo surfaces}
\author[]{Grigory Belousov}
\thanks{The work was partially supported by grant N.Sh.-1987.2008.1}
\maketitle

\begin{abstract}
We prove that a del Pezzo surface with Picard number one has at most
four singular points.
\end{abstract}

\section{intoduction}

A \emph{log del Pezzo surface} is a projective algebraic surface $X$
with only quotient singularities and ample anticanonical divisor
$-K_{X}$.

Del Pezzo surfaces naturally appear in the log minimal model program
(see, e. g., {\cite{KMM}}). The most interesting class of del Pezzo
surfaces is the class of surfaces with Picard number $1$. It is
known that a log del Pezzo surface of Picard number one  has at most
five singular points (see {\cite{KeM}}). Earlier the author proved
there is no log del Pezzo surfaces of Picard number one with five
singular points {\cite{Bel}}. In this paper we give another, simpler
proof.

\begin{theorem}
\label{1.2}
Let $X$ be a del Pezzo surface with log terminal singularities and
Picard number is $1$. Then $X$ has at most four singular points.
\end{theorem}

Recall that a normal complex projective surface is called
a \textit{rational homology
projective plane} if it has the same Betti numbers as the projective
plane $\mathbb P^2$.
J. Koll\'ar {\cite{Kol} posed the problem to find rational homology
$\mathbb P^2$'s with quotient singularities having five singular points.
In  {\cite{HK}} this problem is solved for the case of
numerically effective $K_{X}$.
Our main theorem solves Koll\'ar's problem in the case where $-K_X$ is negative.

The author is grateful to Professor Y. G. Prokhorov for suggesting
me this problem and for his help.

\section{Preliminary results}

We work over complex number field $\CC$.
We employ the following notation:

\begin{itemize}
\item
$(-n)$-curve is a smooth rational curve with self intersection
number $-n$.
\item
$K_{X}$: the canonical divisor on $X$.
\item
$\rho(X)$: the Picard number of $X$.
\end{itemize}

\begin{theorem}[see {\cite[Corollary 9.2]{KeM}}]
\label{2.1} Let $X$ be a rational surface with log terminal
singularities and $\rho(X)=1$. Then
\[
\sum_{P\in X}\frac{m_P-1}{m_P}\leq 3,\leqno{(*)}
\]
where
$m_P$ is the order of the local fundamental group $\pi_1(U_P-\{P\})$
($U_P$ is a sufficiently small neighborhood of $P$).
\end{theorem}

So, every rational surface $X$ with log terminal singularities and
Picard number one has at most six singular points. Assume that $X$
has exactly six singular points. Then by $(*)$ all singularities are
Du Val. This contradicts the classification of del Pezzo surfaces
with Du Val singularities (see, e. g., {\cite{Fur}}, {\cite{Ma}}).

\begin{zero}
Thus to prove Theorem \ref{1.2} it is sufficient to show that there
is no log del Pezzo surfaces with five singular points and Picard
number one. Assume the contrary: there is log del Pezzo surfaces
with five singular points and Picard number one. Let $P_{1},\dots
,P_{5}\in X$ be singular points and $U_{P_{i}}\ni P_{i}$ small
analytic neighborhood. By Theorem \ref{2.1} the collection of
orders of groups $\pi_{1}(U_{P_{1}}-{P_{1}}),\dots
\pi_{1}(U_{P_{5}}-{P_{5}})$ up to permutations is one of the
following:
\begin{case}
$(2, 2, 3, 3, 3)$, $(2, 2, 2, 4, 4)$, $(2, 2, 2, 3, n)$, \ $n=3$, $4$, $5$, $6$,
\end{case}
\begin{case}
$(2, 2, 2, 2, n)$, \quad
$n\geq 2$.
\end{case}
\end{zero}

\begin{remark}
According to the classification of del Pezzo surfaces with Du Val
singularities we may assume that there is a non-Du Val singular
point. The case 2.2.1 is discussed in {\cite[proof of Prop.
6.1]{HK}}. Thus it is sufficient to consider case 2.2.2.
\end{remark}

\begin{zero}
\label{2.3}\textbf{Notation and assumptions.} Let $X$ be a del Pezzo
surface with log terminal singularities and Picard number
$\rho(X)=1$. We assume that we are in case 2.2.2, i. e. the singular
locus of $X$ consists of four points $P_{1}, P_{2}, P_{3}, P_{4}$ of
type $A_{1}$ and one more non Du Val singular point $P_{5}$ with
$|\pi_{1}(U_{P_{5}}-{P_{5}})|=n\geq 3$. Let $\pi\colon \bar{X}
\rightarrow X$ be the minimal resolution and let $D=\sum_{i=1}^{n}
D_{i}$ be the reduced exceptional divisor, where the $D_{i}$ are
irreducible components. Then there exists a uniquely defined an
effective $\QQ$-divisor $D^\sharp=\sum_{i=1}^{n} \alpha_{i}D_{i}$
such that $\pi^*(K_{X})\equiv D^\sharp+K_{\bar{X}}$.
\end{zero}

\begin{lemma}[see, e. g., {\cite[Lemma 1.5]{Za}}]
\label{2.4} Under the condition of \ref{2.3}, let $\Phi\colon \bar
X\rightarrow \PP^1$ be a generically $\PP^1$-fibration. Let $m$ be a
number of irreducible components of $D$ not contain in any fiber of
$\Phi$ and let $d_{f}$ be a number of $(-1)$-curves contained in a
fiber $f$. Then
\begin{enumerate}
\item
$m=1+\sum_{f}(d_{f}-1)$.
\item If $E$ is a unique $(-1)$-curve in $f$, then the coefficient $E$ in $f$ is at least two.
\end{enumerate}
\end{lemma}

The following lemma is a consequence of the Cone Theorem.

\begin{lemma}[see, e. g., {\cite[Lemma 1.3]{Za}}]
\label{2.5} Under the condition of \ref{2.3}, every curve on
$\bar{X}$ with negative selfintersection number is either
$(-1)$-curve or a component of $D$.
\end{lemma}

\begin{definition}
Let $(Y,D)$ be a projective log surface. $(Y,D)$ is called the
\emph{weak log del Pezzo} surface if the pair $(Y,D)$ is klt and the
divisor $-(K_{Y}+D)$ is nef and big.
\end{definition}

For example, in the above notation, $(\bar {X}, D^\sharp)$ is a weak
del Pezzo surface. Note that if $(Y,D)$ is a weak log del Pezzo
surface with $\rho(Y)=1$ then divisor $-(K_{Y}+D)=A$ is ample and
$Y$ has only log terminal singularities. Hence, $Y$ is a log del
Pezzo surface.

\begin{lemma}[see, e. g., {\cite[Lemma 2.9]{Bel}}]
\label{2.7} Suppose $(Y,D)$ is a weak log del Pezzo surface. Let
$f\colon Y\to Y'$ be a birational contraction. Then $(Y',D'=f_{*}D)$
is also a weak log del Pezzo surface.
\end{lemma}

\section{Proof of the main theorem: the case where $X$ has cyclic quotient singularities}
In this section we assume that $X$ has only cyclic quotient
singularities.

The following lemma is very similar to that in {\cite{HP}}. For
convenience of the reader we give a complete proof.

\begin{lemma}
\label{3.1} Under the condition of \ref{2.3}, suppose $P_{5}$ is a
cyclic quotient singularity. Then there exists a generically
$\PP^1$-fibration $\Phi:\bar{X}\rightarrow\PP^1$ such that $f\cdot
D\leq 2$, where $f$ is a fiber of $\Phi$.
\end{lemma}

\begin{proof}
Let $\nu:\hat{X}\rightarrow X$ be the minimal resolution of the non
Du Val singularities and let $E=\sum E_{i}$ be the exceptional
divisor. By {\cite[Corollary 1.3]{PV}} or {\cite[Lemma 10.4]{KeM}}
we have $|-K_{X}|\neq\varnothing$. Take $B\in|-K_{X}|$. Then we can
write
$$K_{\hat{X}}+\hat{B}=\nu^*(K_{X}+B)\sim 0,$$ where $\hat{B}$ is an
effective integral divisor. We obviously have $\hat{B}\geq E$.

Run the MMP on $\hat{X}$. We obtain a birational morphism
$\phi:\hat{X}\rightarrow\tilde{X}$ such that $\tilde{X}$ has only Du
Val singularities and either $\rho(\tilde{X})=2$ and there is a
generically $\PP^1$-fibration $\psi:\tilde{X}\rightarrow\PP^1$ or
$\rho(\tilde{X})=1$. Moreover, $\phi$ is a composition
$$\begin{CD}\hat{X}=X_{1}@>\phi_{1}>>X_{2}@>\phi_{2}>>\dots@>\phi_{n}>>X_{n+1}=\tilde{X},\end{CD}$$
where $\phi_{i}$ is a weighted blowup of a smooth point of $X_{i+1}$
with weights $(1, n_{i})$ (see {\cite{Mor}}).

Assume that $\rho(\tilde{X})=1$, then every singular point on
$\tilde{X}$ is of type $A_{1}$. By the classification of del Pezzo
surfaces with Du Val singularities and Picard number one (see, e.
g., {\cite{Fur}}, {\cite{Ma}}) we have $\tilde{X}=\PP^2$ or
$\tilde{X}=\PP(1,1,2)$.

Assume that $\rho(\tilde{X})=1$ and $\tilde{X}=\PP(1,1,2)$. Since
$\phi_{*}(\hat{B})$ has at most two components, we see that $\phi$
contracts at most two curves $K_{1}$ and $K_{2}$ such that $K_{i}$
is not component of $E$. Since $X$ has four singular points of type
$A_{1}$, we see that $\tilde{X}$ has at least two singular points, a
contradiction.

Assume that $\rho(\tilde{X})=1$ and $\tilde{X}=\PP^2$. Since
$\phi_{*}(\hat{B})$ has at most three components, we see that $\phi$
contracts at most three curves $K_{1}$, $K_{2}$ and $K_{3}$ such
that $K_{i}$ is not component of $E$. Since $X$ has four singular
points of type $A_{1}$, we see that $\tilde{X}$ has at least one
singular point, a contradiction.

Therefore, $\rho(\tilde{X})=2$ and there is a generically
$\PP^1$-fibration $\psi:\tilde{X}\rightarrow\PP^1$. Let
$g:\bar{X}\rightarrow\hat{X}$ be the minimal resolution of
$\hat{X}$. Let $\Phi'=\psi\circ\phi$ and let $f'$ be a fiber of
$\Phi'$. Then $f'\cdot E\leq f'\cdot\hat{B}=-K_{\hat{X}}\cdot f'=2$.
Set $\Phi=\Phi'\circ g$.
\end{proof}

\begin{zero}
Let $f$ be a fiber of $\Phi$. By Lemma \ref{3.1} we have the
following cases:
\begin{case}
$f$ meets exactly one irreducible component $D_{0}$ of $D$ and
$f\cdot D_{0}=1$.

Let $L$ be a singular fiber of $\Phi$. By Lemma \ref{2.4} (1) the
fiber $L$ contains exactly one $(-1)$-curve $F$. By Lemma \ref{2.4}
(2) $F$ does not meet $D_{0}$. Then $F$ meets at most two components
of $D$. Blowup one of the points of intersection $F$ and $D$. We
obtain a surface $Y$. Let $h:Y\rightarrow Y'$ be a contraction of
all curves with selfintersection number at most $-2$. Note that $Y'$
has only log terminal singularities but not of type 2.2.2, a
contradiction.
\end{case}

\begin{case}
$f$ meets exactly two irreducible components $D_{1}$, $D_{2}$ of $D$
and $D_{1}\cdot f=D_{2}\cdot f=1$.

By Lemma \ref{2.4} (1) there exists a unique singular fiber $L$ such
that $L$ has two $(-1)$-curves $F_{1}$ and $F_{2}$. Note that one of
this curves, say $F_{1}$, meets $D$ at one or two points. Blowup one
of the points of intersection $F_{1}$ and $D$. We obtain a surface
$Y$. Let $h:Y\rightarrow Y'$ be a contraction of all curves with
selfintersection number at most $-2$. Note that $Y'$ has only log
terminal singularities but not of type 2.2.2, a contradiction.
\end{case}

\begin{case}
$f$ meets exactly one irreducible component $D_{0}$ of $D$ and
$f\cdot D_{0}=2$. Let $A$ be a connected component of $D$ containing
$D_{0}$.

By Lemma \ref{2.4} (1) every singular fiber of $\Phi$ contains
exactly one $(-1)$-curve. Note that every singular fiber of $\Phi$
either contains two connected components of $A-D_{0}$ or the
coefficient of a unique $(-1)$-curve in this fiber is equal two. If
a singular fiber $L$ contains exactly one $(-1)$-curve with
coefficient two, then the dual graph of $L$ is following:
\begin{equation*}
 \xymatrix
@R=1pc @C=1.8pc { \overset{-2}\circ\ar@{-}[r]&
\overset{-1}\circ\ar@{-}[r]& \overset{-2}\circ } \leqno{(**)}
\end{equation*}

Since $X$ has five singular points with orders of local fundamental
groups $(2,2,2,2,n)$, we see that $\Phi$ has two singular fibers
$L_{1}$, $L_{2}$ of type $(**)$ and possibly one more singular fiber
$L_{3}$. Note that $L_{3}$ contains both connected component of
$A-D_{0}$. Let $\mu:\bar{X}\rightarrow\FF_{n}$ be the contraction of
all $(-1)$-curves in fibers of $\Phi$, where $\FF_{n}$ is the
Hirzebruch surface of degree $n$ (rational ruled surface) and
$n=0,1$. Denote $\tilde{D_{0}}:=\mu_{*}D_{0}$. Note that
$\tilde{D_{0}}\sim 2M+kf$, where $M^2=-n$ and $M\cdot f=1$. Since we
contract at most five curves that meet $D_{0}$, and $D_{0}^2\leq
-2$, we see that $0<\tilde{D_{0}}^2\leq 3$. Hence, $0<-4n+4k\leq 3$.
This is impossible, a contradiction.
\end{case}
\end{zero}

\section{Proof of the main theorem: the case where $X$ has non-cyclic quotient singularity}

Under the condition of \ref{2.3}, assume $X$ has a non-cyclic
singular point, say $P$. Then there is a unique component $D_{0}$ of
$D$ such that $D_{0}\cdot(D-D_{0})=3$ (see \cite{Br}).

\begin{lemma} There is a generically $\PP^1$-fibration
$\Phi:\bar{X}\rightarrow\PP^1$ such that $\Phi$ has a unique section
$D_{0}$ in $D$ and $D_{0}\cdot f\leq 3$, where $f$ is a fiber of
$\Phi$.
\end{lemma}

\begin{proof}
Recall that $P$ is not Du Val. Let $h:\bar{X}\rightarrow\hat{X}$ be
contract all curves in $D$ except $D_{0}$. Let
$\hat{D_{0}}=h_{*}(D_{0})$ then $\hat{X}$ has seven singular points,
$\rho(\hat{X})=2$ and there is $\nu: \hat{X}\rightarrow X$ such that
$K_{\hat{X}}+a\hat{D_{0}}=\nu^*K_{X}$. Note that $(\hat{X}, aD_{0})$
is a weak log del Pezzo. Let $R$ be the extremal rational curve
different from $\hat{D}$. Let $\phi: \hat{X}\rightarrow\tilde{X}$ be
the contraction of $R$.

\begin{zero}
There are two cases:

\begin{case}
$\rho(\tilde{X})=1$. Then, by Lemma \ref{2.7}, $\tilde{X}$ is a del
Pezzo surface. If the number of singular points of $\hat{X}$ on $R$
is at most two, $\tilde{X}$ has at least five singular points and
all points are cyclic quotients. Thus assume that there is at least
three singular points of $\hat{X}$ on $R$, say $P_{1}$, $P_{2}$,
$P_{3}$. Let $R_{1}=\sum_{i}R_{1i}$, $R_{2}=\sum_{i}R_{2i}$ and
$R_{3}=\sum_{i}R_{3i}$ be the exceptional divisors on $\bar{X}$ over
$P_{1}$, $P_{2}$ and $P_{3}$, respectively. Let $\bar{R}$ is the
proper transformation  of $R$ on $\bar{X}$. Since $\bar{R}$ is not
component of $D$, we see that $\bar{R}^2\geq-1$. Indeed, this
follows from Lemma \ref{2.5}. Note that matrix of intersection of
component $\bar{R}+R_{1}+R_{2}+R_{3}$ is not negative definite.
Hence, $\bar{R}+E_{1}+E_{2}+E_{3}$ can not be contracted, a
contradiction.
\end{case}

\begin{case}
$\tilde{X}=\PP^1$. Let $g:\bar{X}\rightarrow\hat{X}$ be the
resolution of singularities. Then $\Phi=\phi\circ g:
\bar{X}\rightarrow\PP^1$. Note that there is a unique horizontal
curve $D_{0}$ in $D$. Let $f$ be a fiber of $\Phi$. Denote
coefficient of $D_{0}$ in $D^\sharp$ by $\alpha$. Then
\[
0>(K_{\bar{X}}+D^\sharp)\cdot f=-2+\alpha (D_{0}\cdot f).
\]
Hence, $D_{0}\cdot f<\frac{2}{\alpha}$. Since $P$ is not Du Val, we
see that $\alpha\geq\frac{1}{2}$. Hence, $D_{0}\cdot f\leq 3$.
\end{case}
\end{zero}
\end{proof}

By Lemma \ref{2.4} (1) every singular fiber of $\Phi$ contains
exactly one $(-1)$-curve. Let $B$ be the exceptional divisor
corresponding to the non-cyclic singular point. Note that $B$
contains $D_{0}$.

\begin{zero}
Consider three cases.
\begin{case}
$D_{0}\cdot f=1$. Then every singular fiber of $\Phi$ contains
exactly one connected component of $B-D_{0}$. On the other hand,
$B-D_{0}$ contains three connected component. Hence $X$ has at most
four singular points, a contradiction.
\end{case}
\begin{case}
$D_{0}\cdot f=2$. Let $F_{1}$, $F_{2}$, $F_{3}$ be a connected
components of $B-D_{0}$. We may assume $F_{1}$ is $(-2)$-curve (see
\cite{Br}). Let $L_{1}$ be a singular fiber of $\Phi$. Assume that
$L_{1}$ contains $F_{1}$. Then $L_{1}$ is type $(**)$ and $L_{1}$
contain $F_{2}$. Hence, $F_{2}$ is a $(-2)$-curve. Let $L_{2}$ be a
singular fiber of $\Phi$. Assume that $L_{2}$ contains $F_{3}$ and
let $E$ be a unique $(-1)$-curve in $L_{2}$. By blowing up the point
of intersection $E$ and $F_{3}$, we obtain a surface $Y$. Let
$h:Y\rightarrow Y'$ be a contraction of all curve with
selfintersection number at most $-2$. Note that $Y'$ has only log
terminal singularities but not of type 2.2.2, a contradiction.
\end{case}
\begin{case}
$D_{0}\cdot f=3$. Since every component of $D-B$ is a $(-2)$-curve,
we see that every singular fiber of $\Phi$ contains a connected
component of $B-D_{0}$. Note that $B-D_{0}$ contains three connected
components. Hence $X$ has at most four singular points, a
contradiction.
\end{case}
\end{zero}
This completes the proof of Theorem \ref{1.2}.


\begin{thebibliography}{lll}

\bibitem{Bel}
G. N. Belousov, \emph{Del Pezzo surfaces with log terminal
singularities}, Mat. Zametki \textbf{83} (2008), no. 2, 170-180
(Russian) English translate Math. Notes.

\bibitem{Br}
Brieskorn E. \emph{Rationale Singularit\"{a}ten komplexer
Fl\"{a}chen}, Invent. Math. \textbf{4} (1968), 336 -- 358.

\bibitem{Fur}
Furushima M. \emph{Singular del Pezzo surfaces and analytic
compactifications of 3-dimensional complex affine space $C^3$},
Nagoya Math. J. \textbf{104} (1986), 1 -- 28.

\bibitem{HK} D. Hwang, J. Keum \emph{The maximum number of singular points on rational homology projective
planes} arXiv:math.AG/0801.3021v3.

\bibitem{HP} P. Hacking Yu. Prokhorov \emph{Smoothable del Pezzo
surfaces with quotient singularities} Unpublished manuscript.

\bibitem{K}
Kawamata Y. \emph{Crepant blowing-up of 3-dimensional canonical
singularities and its application to degenerations of surfaces},
Ann. of Math. \textbf{127} (1988), 93 -- 163.

\bibitem{KMM}
Kawamata Y., Matsuda K. $\&$ Matsuki J. \emph{Introduction to the
minimal model program}, Adv. Stud. Pure Math. \textbf{10} (1987),
283 -- 360.

\bibitem{KeM}
Keel S. $\&$ McKernan J. \emph{Rational curves on quasi-projective
surfaces}, Memoirs AMS \textbf{140} (1999), no. 669.

\bibitem{Kol}
Kollar J. \emph{Is there a topological Bogomolov-Miyaoka-Yau
inequality?}, Pure and Applied Math. Quarterly \textbf{4} No. 2
(2008).

\bibitem{Ma}
Miyanishi M. $\&$ Zhang D. -Q. \emph{Gorenstein log del Pezzo surfaces of rank one}, J. Algebra. \textbf{118} (1988), 63 -- 84.

\bibitem{Mor}
Morrison D. \emph{The Birational Geometry of Surfaces with Rational
Double Points} Math. Ann. \textbf{271} (1985), 415-438.
\bibitem{PV}
Prokhorov Yu. G. Verevkin A. B. \emph{The Riemann–Roch theorem on
surfaces with log terminal singularities} J. Math Sci. (N. Y.)
\textbf{140} (2007), no. 2, 200-205.
\bibitem{Za}
Zhang D.-Q. \emph{Logarithmic del Pezzo surfaces of rang one with
contractible boundaries}, Osaka J. Math. \textbf{25} (1988), 461 --
497.

\end{thebibliography}
\end{document}